\documentclass[vecphys]{svmult}
\usepackage{graphicx}        
\usepackage{multicol}        
\usepackage[bottom]{footmisc}
\usepackage{amsmath,amsfonts}

\begin{document}



\DeclareGraphicsExtensions {.eps}

\title*{Solving One Dimensional Scalar Conservation Laws by Particle Management}
\titlerunning{Solving Conservation Laws by Particle Management}
\author{Yossi Farjoun\and Benjamin Seibold}
\institute{Department of Mathematics,
Massachusetts Institute of Technology,
77~Massachusetts Avenue,
Cambridge MA 02139, USA \\
\texttt{\{yfarjoun,seibold\}@math.mit.edu}}
\maketitle
\begin{abstract}
We present a meshfree numerical solver for scalar conservation laws in one space
dimension. Points representing the solution are moved according to their characteristic
velocities.  Particle interaction is resolved by purely local particle management.
Since no global remeshing is required, shocks stay sharp and propagate at the correct
speed, while rarefaction waves are created where appropriate. The method is TVD, entropy
decreasing, exactly conservative, and has no numerical dissipation. Difficulties
involving transonic points do not occur, however inflection points of the flux function
pose a slight challenge, which can be overcome by a special treatment. Away from shocks
the method is second order accurate, while shocks are resolved with first order accuracy.
A postprocessing step can recover the second order accuracy. The method is compared to
CLAWPACK in test cases and is found to yield an increase in accuracy for comparable
resolutions.
\end{abstract}
\begin{keywords}
conservation law, meshfree, particle management, entropy
\end{keywords}

\newcommand{\prn}[1]{\left(#1\right)}
\newcommand{\brk}[1]{\left[#1\right]}
\newcommand{\abs}[1]{\left|#1\right|}
\newcommand{\ud}[1]{\, \mathrm{d}#1}
\newcommand{\pd}[2]{\frac{\partial#1}{\partial#2}}
\newcommand{\sign}{{\rm sign}}
\renewcommand{\smartqed}{\hspace*{\fill}\qed}

\section{Introduction}
Lagrangian particle methods approximate the solution of flow equations using a cloud of
points which move with the flow. Examples are vortex methods \cite{seibold::Chorin1973},
smoothed particle hydrodynamics (SPH) \cite{seibold::Lucy1977,seibold::GingoldMonaghan1977},
or generalized SPH methods \cite{seibold::Dilts1999}. The latter are typically based
on generalized meshfree finite difference schemes \cite{seibold::LiszkaOrkisz1980}.
An example is the finite pointset method (FPM) \cite{seibold::KuhnertTiwari2002}.
Moving the computational nodes with the flow velocity $\vec{v}$ allows the discretization
of the governing equations in their more natural Lagrangian frame of reference. The
material derivative \mbox{$\frac{D}{Dt} = \pd{}{t}+\vec{v}\cdot\nabla$} becomes a simple
time derivative. For a conservation law, the natural velocity is the characteristic
velocity. In a frame of reference which is moving with this velocity, the equation
states that the function value remains constant. Of course, this is only valid where the
solution is smooth. In this case, characteristic particle methods are very simple
solution methods for conservation laws. 

In spite of the obvious advantages of particle methods, almost all numerical methods for
conservation laws operate on a fixed Eulerian grid, even though significant work has to
be invested to solve even a simple advection problem preserving  sharp features and
without creating oscillations. Leaving aspects of implementation complexity aside, two
main reasons favor fixed grid methods: First, a fixed grid allows an easy generalization
to higher space dimensions using dimensional splitting. Second, in particle methods one
has to deal with the interaction of characteristics. While the former point remains
admittedly present, the latter aspect is addressed in this contribution.

Most methods which use the characteristic nature of the conservation law circumvent
the problem of crossing characteristics by remeshing. Before any particles can interact,
the numerical solution is interpolated onto ``nicely'' distributed particles, for
instance onto an equidistant grid -- in which case the approach is essentially a fixed
grid method. The CIR-method \cite{seibold::CourantIsaacsonRees1952} is an example.
Such approaches incur multiple drawbacks: First, the shortest interaction time determines
the global time step. Second, the error due to the global interpolation may yield
numerical dissipation and dispersion. Finally, such schemes are not conservative when
shocks are present. In practice, finite volume approaches, such as
Godunov methods with appropriate limiters \cite{seibold::VanLeer1974}, or
ENO \cite{seibold::HartenEngquistOsherChakravarthy1987}/WENO
\cite{seibold::LiuOsherChan1994} schemes are used to compute weak entropy solutions
that show neither too much oscillations nor too much numerical dissipation.

With moving particles, two fundamental problems arise: On the one hand, neighboring
particles may depart from each other, resulting in poorly resolved regions. On the other
hand, a particle may (if left unchecked) overtake a neighbor, which results in a
``breaking wave'' solution. The former problem can be remedied by inserting particles. 
The latter has to be resolved by merging particles. When characteristic particles
interact (i.e.~one overtakes the other) one is dealing with a shock, thus particles must
be merged.

In this contribution, we present a local and conservative particle management (inserting
and merging particles) that yields no numerical dissipation (where solutions are smooth)
and correct shock speeds (where they are not). The particle management is based on exact
conservation properties between neighboring particles, which are derived in
Sect.~\ref{seibold::sec:conservation_laws}.
In Sect.~\ref{seibold::sec:method_description}, we outline our numerical method.
The heart of our method, the particle management, is derived in
Sect.~\ref{seibold::sec:interp_particle_management}.
There, we also show that the method is TVD.
In Sect.~\ref{seibold::sec:entropy}, we prove that the numerical
solutions satisfy the Kru\v zkov entropy condition, thus showing that
the solutions we find are entropy solutions for any convex entropy function.
In Sect.~\ref{seibold::sec:numerical_results}, we apply the method to examples and
compare it to traditional finite volume methods using CLAWPACK.
In Sect.~\ref{seibold::sec:inflection_points} we present how
non-convex flux functions can be treated.
Finally, in Sect.~\ref{seibold::sec:outlook}, we outline possible extensions
and conclusions. These include applications of the 1D solver itself as well as possible
extensions beyond the 1D scalar case.

\section{Scalar Conservation Laws}
\label{seibold::sec:conservation_laws}
Consider a one-dimensional scalar conservation law
\begin{equation}
u_t+\prn{f(u)}_x = 0, \quad u(x,0) = u_0(x)
\label{seibold::eq:conslaw}
\end{equation}
with $f'$ continuous.
As long as the solution is smooth, it can be obtained by the method of
characteristics \cite{seibold::Evans1998}.
The function $u(x(t), t)$ is constant along the characteristic curve
\begin{equation}
x(t) = x(0)+f'(u(x(0),0))\,t \;.
\label{seibold::eq:characteristic}
\end{equation}
For nonlinear functions $f$ the characteristic curves can ``collide'', resulting in a
shock, whose speed is given by the Rankine-Hugoniot condition \cite{seibold::Evans1998}.
Discontinuities are shocks only if the characteristic curves run into them. 
Other discontinuities become rarefaction waves, i.e.~continuous functions which attain
every value between the left and the right states. If the flux function $f$ is convex or
concave between the left and right state of a discontinuity, then the solution is either
a shock or a rarefaction. If $f''$ switches sign between between the two states, then
a combination of a shock and a rarefaction occur. These physical solutions are defined
by a weak formulation of \eqref{seibold::eq:conslaw} accompanied by an entropy condition.

\subsection{Conservation Properties}
\label{seibold::subsec:cons_properties}
Conservations laws conserve the total area under the solution
\begin{equation}
\frac{d}{dt}\int\limits_{-\infty}^\infty u(x,t)\ud{x} = 0 \;.
\label{seibold::eq:cons_law_area}
\end{equation}
The change of area between two \emph{moving} points $b_1(t)$ and $b_2(t)$ is given by
\begin{align*}
\frac{d}{dt}\int_{b_1(t)}^{b_2(t)} u(x,t)\ud{x}
= b_2'(t)\,u(b_2(t),t)-b_1'(t)\,u(b_1(t),t)+\int_{b_1(t)}^{b_2(t)} u_t(x,t)\ud{x} \\
= \prn{b_2'(t)\,u(b_2(t),t)-b_1'(t)\,u(b_1(t),t)}-\prn{f(u(b_2(t),t))-f(u(b_1(t),t))} \;.
\end{align*}
If $x_1(t)$ and $x_2(t)$ are \emph{characteristic} points,
that is, points following the characteristics of a smooth solution as
in equation \eqref{seibold::eq:characteristic}, we have that
$x_1'(t)=f'(u_1)$ and $x_2'(t)=f'(u_2)$.
Therefore, the change of area between $x_1$ and $x_2$ is
\begin{equation}
\prn{f'(u_2)u_2-f'(u_1)u_1}-\prn{f(u_2)-f(u_1)} = \brk{f'(u)u-f(u)}_{u_1}^{u_2} \;,
\label{seibold::eq:area_deriv}
\end{equation}
where $\brk{g(x)}_a^b=g(b)-g(a)$.
Equation \eqref{seibold::eq:area_deriv} implies that the change of area between two
characteristic points does \emph{not} depend on the positions of the points, only on the
left state $u_1$ and right state $u_2$ and the flux function. Since the two states do
not change as the points move, the area between the two points changes linearly, as does
the distance between them:
\begin{equation}
\frac{d}{dt}(x_2(t)-x_1(t)) = x_2'(t)-x_1'(t) = f'(u_2)-f'(u_1)
= \brk{f'(u)}_{u_1}^{u_2} \;.
\label{seibold::eq:distance_deriv}
\end{equation}
If the two points move at different speeds, then there is a time $t_0$ (which may be
larger or smaller than $t$) at which they have the same position. 
Thus at time $t=t_0$ the distance between them, and the
area between them equal zero. 
From \eqref{seibold::eq:area_deriv} and \eqref{seibold::eq:distance_deriv} we have that
\begin{align*}
\int_{x_1(t)}^{x_2(t)} u(x,t)\ud{x} &= (t-t_0) \cdot \brk{f'(u)u-f(u)}_{u_1}^{u_2} \;, \\
x_2(t)-x_1(t) &= (t-t_0) \cdot \brk{f'(u)}_{u_1}^{u_2} \;.
\end{align*}
In short, the area between two Lagrangian points can be written as
\begin{equation}
\int_{x_1(t)}^{x_2(t)} u(x,t)\ud{x} = (x_2(t)-x_1(t))\,a(u_1,u_2) \;,
\end{equation}
where $a(u_1, u_2)$ is the nonlinear average function
\begin{equation}
a(u_1,u_2)
= \frac{\brk{f'(u)u-f(u)}_{u_1}^{u_2}}{\brk{f'(u)}_{u_1}^{u_2}}
= \frac{\int_{u_1}^{u_2}f''(u)\,u \ud{u}}{\int_{u_1}^{u_2}f''(u) \ud{u}} \;.
\label{seibold::eq:nonlinear_average}
\end{equation}
The integral form shows that $a$ is indeed an average of $u$, weighted by $f''$.
This last observation needs one additional assumption: that the points $x_1$
and $x_2$ remain characteristic point between $t$ and $t_0$. That is,
that a shock does not develop between the two points before $t_0$.
Our numerical method relies heavily on the nonlinear average $a(\cdot,\cdot)$.
\begin{lemma}
\label{seibold::thm:average_properties}
Let $f$ be strictly convex or concave in $[u_L,u_U]$, that is
$f''<0$ or $f''>0$ in $(u_L,u_U)$.
Then for all $u_1, u_2\in [u_L,u_U]$, the average \eqref{seibold::eq:nonlinear_average}
is\dots
\begin{enumerate}
\item the same for $f$ and $-f$. Thus we assume $f''>0$ WLOG;
\item symmetric, $a(u_1,u_2) = a(u_2,u_1)$. Thus we assume $u_1\le u_2$ WLOG;
\item an average, i.e.~$a(u_1,u_2)\in(u_1,u_2)$, for $u_1\neq u_2$;
\item strictly increasing in both $u_1$ and $u_2$; and
\item continuous at $u_1=u_2$, with $a(u_1,u_1)=u_1$.
\end{enumerate}
\end{lemma}
\noindent
\begin{proof}
We only include here the proof of 4.
We show that $a(u_1,u_2)$ is strictly increasing in the second argument.
Let $u_1<u_2<u_3$, $u_i\in [u_L, u_U]$. Then
\begin{align*}
a(u_1, u_3)
&=\frac{\int_{u_1}^{u_2}f''(u)u\ud{u}+\int_{u_2}^{u_3}f''(u)u\ud{u}}
{\int_{u_1}^{u_3}f''(u)\ud{u}} \\
&>\frac{a(u_1,u_2)\int_{u_1}^{u_2}f''(u)\ud{u}+a(u_1,u_2)\int_{u_2}^{u_3}f''(u)\ud{u}}
{\int_{u_1}^{u_3}f''(u)\ud{u}}
= a(u_1, u_2) \;.
\end{align*}
Similar arguments show the result for the first argument.\smartqed
\end{proof}

\section{Description of the Particle Method}
\label{seibold::sec:method_description}
The first step is to approximate the initial function $u_0$ by a finite number of points
$x_1,\dots,x_m$ with function values $u_1,\dots,u_m$. A straightforward strategy is to
place $x_1,\dots,x_m$ equidistant on the interval of interest and assign $u_i = u_0(x_i)$.
More efficient adaptive sampling strategies can be used, since our method does not impose
any requirements on the point distribution. For instance, one can choose $x_i$ and $u_i$
to minimize the $L^1$ error, using the specific interpolation introduced in
Sect.~\ref{seibold::sec:interp_particle_management}. This strategy is the topic of future
work. The points are ordered so that $x_1<\dots<x_m$.
The evolution of the solution is found by moving each point $x_i$ with speed $f'(u_i)$.
This is possible as long as there are no ``collisions'' between points. Two neighboring
points $x_i(t)$ and $x_{i+1}(t)$ collide at time \mbox{$t+\varDelta t_i$}, where
\begin{equation}
\varDelta t_i = -\frac{x_{i+1}-x_i}{f'(u_{i+1})-f'(u_i)} \;.
\label{seibold::eq:intersection_time}
\end{equation}
A positive $\varDelta t_i$ indicates that the two points will eventually collide.
Thus, $t+\varDelta t_{\text s}$ is the time of the next particle
collision\footnote{If the set $\{i|\varDelta t_i\ge 0\}$ is empty,
then $\varDelta t_{\text s}=\infty$.}, where
\begin{equation*}
\varDelta t_{\text s} = \min_i\{\varDelta t_i | \varDelta t_i\ge 0\} \;.
\end{equation*}
For any time increment $\varDelta t\le \varDelta t_{\text s}$ the points can be moved
directly to their new positions $x_i+f'(u_i)\varDelta t$.  Thus, we can step forward an
amount $\varDelta t_{\text s}$, and move all points accordingly. Then, at least one
particle will share its position with another. To proceed further, we merge each such
pair of particles.
If the collision time $\varDelta t_i$ is negative, the points depart from each other.
Although at each of the points the correct function value is preserved, after some
time their distance may be unsatisfyingly large, as the amount of error
introduced during a merge grows with the size of the gaps in the
neighboring particles. To avoid this, we insert new points into large
gaps between points \emph{before} merging particles.
In Sect.~\ref{seibold::subsec:particle_management} we derive positions
and values of the new particles that assure that the method is
conservative, TVD, and entropy diminishing.

\section{Interpolation and Particle Management}
\label{seibold::sec:interp_particle_management}
The movement of the particles is given by a fundamental property of the conservation
law \eqref{seibold::eq:conslaw}: its characteristic equation
\eqref{seibold::eq:characteristic}. We derive particle management to satisfy another
fundamental property: the conservation of area \eqref{seibold::eq:cons_law_area}.
Using the conservation principles derived in
Sect.~\ref{seibold::sec:conservation_laws}, the function value of an inserted
or merged particle is chosen, such that area is conserved exactly. 
A simple condition on the particles guarantees that the entropy does not increase.
In addition, we define an interpolating function between two neighboring particles,
so that the change of area satisfies relation \eqref{seibold::eq:area_deriv}.
Furthermore, this interpolation is an analytical solution to the conservation law.

\subsection{Conservative Particle Management}
\label{seibold::subsec:particle_management}
\begin{figure}
\centering
\begin{minipage}[t]{.47\textwidth}
\centering
\includegraphics[width=0.99\textwidth]{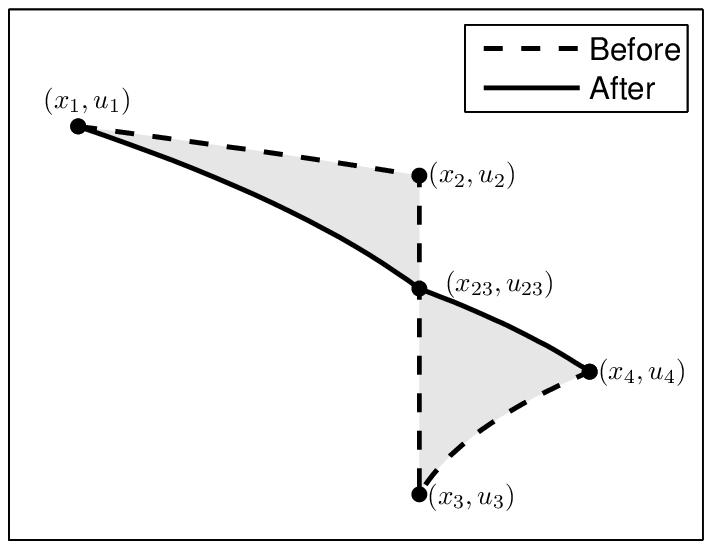}
\caption{Merging two particles}
\label{seibold::fig:merging}
\end{minipage}
\hfill
\begin{minipage}[t]{.495\textwidth}
\centering
\includegraphics[width=0.99\textwidth]{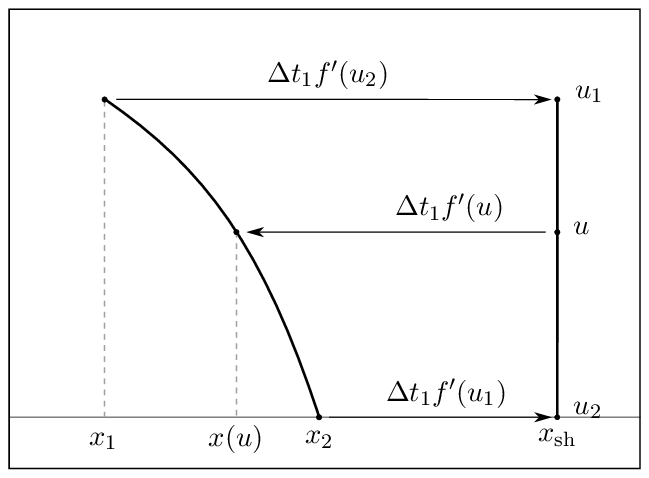}
\caption{Definition of the interpolation}
\label{seibold::fig:def_interp}
\end{minipage}
\end{figure}

Consider four neighboring particles located at \mbox{$x_1<x_2\le x_3<x_4$} with
associated function values $u_1$, $u_2$, $u_3$, $u_4$. Assume that the flux $f$ is
strictly convex or concave on the range of function values $[\min_i(u_i),\max_i(u_i)]$.
If $u_2\neq u_3$, the particles' velocities must differ $f'(u_2)\neq f'(u_3)$, which
gives rise to two possible cases that require particle management:
\begin{itemize}
\item \textbf{Inserting:}
The two particles deviate, i.e.~$f'(u_2)<f'(u_3)$. If the distance $x_3-x_2$ is
larger than a predefined maximum distance $d_{\text max}$, we insert a new particle
$(x_{23},u_{23})$ with $x_2<x_{23}<x_3$ and $u_{23}$ chosen so that the area between
$x_2$ and $x_3$  is preserved by the insertion:
\begin{equation}
(x_{23}-x_2)\,a(u_2,u_{23})+(x_3-x_{23})\,a(u_{23},u_3) = (x_3-x_2)\,a(u_2,u_3) \;.
\label{seibold::eq:area_cond_insert}
\end{equation}
This condition defines a function, connecting $(x_2,u_2)$ with $(x_3,u_3)$, on which
the new particle has to lie. This function is the interpolation defined in
Sect.~\ref{seibold::subsec:interpolation} and illustrated in
Fig.~\ref{seibold::fig:def_interp}.
\item \textbf{Merging:}
The two particles collide, i.e.~$f'(u_2)>f'(u_3)$. If the distance $x_3-x_2$ is
smaller than a preset value $d_{\text min}$ ($d_{\text min}=0$ is possible), we replace
both with a new particle $(x_{23},u_{23})$. The position of the new particle $x_{23}$
is chosen with \mbox{$x_2<x_{23}<x_3$} and $u_{23}$ is chosen so that the total area
between $x_1$ and $x_4$ is unchanged:
\begin{align}
(x_{23}-x_1)\,a(u_1,u_{23})+(x_4-x_{23})\,a(u_{23},u_4)& \nonumber \\
= (x_2\!-\!x_1)\,a(u_1,u_2)+(x_3\!-\!x_2)\,a(u_2,u_3)&+(x_4\!-\!x_3)\,a(u_3,u_4) \;.
\label{seibold::eq:area_cond_merge}
\end{align}
Any particle $(x_{23},u_{23})$ with \mbox{$x_2<x_{23}<x_3$} that satisfies
\eqref{seibold::eq:area_cond_merge} would be a valid choice. 
We choose $x_{23}=\frac{x_2+x_3}{2}$, and then obtain
$u_{23}$ such that \eqref{seibold::eq:area_cond_merge} is satisfied.
Figure~\ref{seibold::fig:merging} illustrates the merging step.
\end{itemize}
Observe that inserting and merging are similar in nature. Conditions
\eqref{seibold::eq:area_cond_insert} and
\eqref{seibold::eq:area_cond_merge} for $u_{23}$  are
nonlinear (unless $f$ is quadratic, see Remark~\ref{seibold::remark:quadratic_flux}).
For most cases $u_{23}=\frac{u_2+u_3}{2}$ is a good initial guess, and
the correct value can be obtained by a few Newton iteration steps. 
The next few claims attest that we can find a unique value $u_{23}$ that satisfies 
\eqref{seibold::eq:area_cond_insert} and \eqref{seibold::eq:area_cond_merge}.
\begin{lemma}
The function value $u_{23}$ for the particle at $x_{23}$ is unique.
\end{lemma}
\begin{proof}
We show the case for merging. The argument for insertion is similar. 
From Lemma~\ref{seibold::thm:average_properties} we have that both $a(u_1,\cdot)$
and $a(\cdot, u_4)$ are strictly increasing. Thus, the LHS of
\eqref{seibold::eq:area_cond_merge} is strictly increasing, and cannot
be the same for different values of $u_{23}$.\smartqed
\end{proof}
\begin{lemma}
\label{seibold::lem:existence_u23}
If $x_2 = x_3 = x_{23}$, there exists $u_{23}\in [u_2,u_3]$ which
satisfies \eqref{seibold::eq:area_cond_merge}. 
\end{lemma} 
\begin{proof}
WLOG we assume that $u_2\le u_3$.
First, we define
\begin{align*}
A\phantom{(u)}&=(x_2-x_1)\,a(u_1,u_2)+(x_4-x_2)\,a(u_3,u_4),\quad \text{and} \\
B(u) &= (x_2-x_1)\,a(u_1,u)+(x_4-x_2)\,a(u,u_4)\;.
\end{align*}
Equation \eqref{seibold::eq:area_cond_merge} is now simply $B(u_{23})=A$.
The monotonicity of $a$ implies that 
\begin{equation}
B(u_2)\le A \le B(u_3) \;.
\end{equation}
Since $a$ is continuous, so is $B$, and the existence of $u_{23}$
follows the intermediate value theorem.\smartqed
\end{proof}
\begin{corollary}
If particles are merged only according to
Lemma~\ref{seibold::lem:existence_u23}, then the total variation of
the solution is either the same as before the merge, or smaller.
\end{corollary}
Merging points only when $x_2=x_3$ can be too restrictive.
Fortunately, the following claim allows for a little more freedom.
\begin{theorem}
\label{seibold::thm:merging_TVD}
Consider four consecutive particles $(x_i,u_i) \ \forall i=1,\ldots,4$.
Merging particles 2 and 3 so that $x_{23}=\frac{x_2+x_3}{2}$ yields
$u_{23}\in[u_2,u_3]$ if 
\begin{equation}
\frac{x_3-x_2}{\abs{u_3-u_2}}
\le\frac{1}{16}\prn{\frac{\min\abs{f''}}{\max\abs{f''}}}^6
\frac{\min\prn{x_4-x_2,x_3-x_1}}{\abs{\max(u_3,u_2)-\min(u_4, u_1)}} \;.
\end{equation}
\end{theorem}
Here the $\min$ and $\max$ of $f''$ are taken over the maximum range of
$u_1,\dots,u_4$. This condition is naturally satisfied if $x_2=x_3$.
\begin{proof}[outline]
The full proof will be given in a future paper. The idea is to merge in two steps:
First, we find a value $\tilde u$ such that setting $u_2=u_3=\tilde u$ (while leaving
$x_2$ and $x_3$ unchanged) results in the same area. Then, we merge the two points
to $u_{23}$. In the first step we bound $\tilde u$ \emph{away} from $u_2$ and $u_3$
(but inside $[u_2,u_3]$), and in the second step we bound $\abs{u_{23}-\tilde u}$
from above.\smartqed
\end{proof}
\begin{theorem}
\label{seibold::thm:arbitrary_times}
The particle method can run to arbitrary times.
\end{theorem}
\begin{proof}
Let $u_L=\min_i u_i$, $u_U=\max_i u_i$, and $C=\max_{[u_L,u_U]}|f''(u)|\cdot(u_U-u_L)$.
For any two particles, one has $|f'(u_{i+1})-f'(u_i)|\le C$. Define
$\varDelta x_i=x_{i+1}-x_i$. After each particle management, the next time increment
(as defined in Sect.~\ref{seibold::sec:method_description}) is at least
$\varDelta t_{\text s}\ge \frac{\min_i{\varDelta x_i}}{C}$. If we do not insert
particles, then in each merge one particle is removed. Hence, a time evolution beyond
any given time is possible, since the increments $\varDelta t_{\text s}$ will increase
eventually. When a particle is inserted (whenever two points are at a distance more
than $d_{\text max}$), the created distances are at least $\frac{d_{\text max}}{2}$,
preserving a lower bound on the following time increment.\smartqed
\end{proof}

\subsection{Conservative Interpolation}
\label{seibold::subsec:interpolation}
The particle management does not require an interpolation between points. As it stands,
it complements the characteristic movement to yield a full particle method for the
conservation law \eqref{seibold::eq:conslaw} that can run for arbitrarily long times.
Yet, for plotting the solution and interpreting approximation properties, it is
desirable to define an interpolation that is compatible with the conservation principles
of the underlying partial differential equation. We define such an interpolation between
each two neighboring points $(x_1,u_1)$ and $(x_2,u_2)$.

In the case $u_1=u_2$, we define the interpolation to be a constant. In the
following, we describe the case $u_1\neq u_2$. Assume that $f$ is strictly convex
or concave in $[u_1,u_2]$. Therefore $f'(u_1)\neq f'(u_2)$. Hence, as derived in
Sect.~\ref{seibold::subsec:cons_properties}, the solution either came from a
discontinuity (i.e.~it is a rarefaction) or it will become a shock. The time
$\varDelta t_1$ until this discontinuity happens is given by
\eqref{seibold::eq:intersection_time}. At time $t+\varDelta t_1$ the points have the
same position $x_1=x_2=x_{\text{sh}}$, as shown in Fig.~\ref{seibold::fig:def_interp}.
At this time the interpolation must be a straight line connecting the two points,
representing a discontinuity at $x_{\text{sh}}$. We require any point of the
interpolating function $(x,u)$ to move with its characteristic velocity $f'(u)$ in
the time between $t$ and $t+\varDelta t_1$. 
This defines the interpolation uniquely as
\begin{equation}
x(u) = x_1-t_1\prn{f'(u)-f'(u_1)}
= x_1+\frac{f'(u)-f'(u_1)}{f'(u_2)-f'(u_1)}(x_2-x_1) \;.
\label{seibold::eq:interpolation_function}
\end{equation}
Defining $x$ as a function of $u$ is in fact an advantage, since at
a discontinuity characteristics for all intermediate values $u$ are defined. Thus,
rarefaction fans arise naturally if $f'(u_1)<f'(u_2)$. Since $f''$ has no inflection
points between $u_1$ and $u_2$, the inverse function $u(x)$ exists. However, it is
only required at a single point for particle management. For plotting purposes we can
always plot $x(u)$ instead.
\begin{lemma}
The interpolation \eqref{seibold::eq:interpolation_function}
is a solution of the conservation law \eqref{seibold::eq:conslaw}.
\end{lemma}
\begin{proof}
Using that \mbox{$\dot x_i(t)=f'(u_i)$} for~$i=1,2$ one obtains
\begin{align*}
\frac{\partial x}{\partial t}(u,t)&= \dot x_1
+\frac{f'(u)-f'(u_1)}{f'(u_2)-f'(u_1)}(\dot x_2-\dot x_1)\\
&=f'(u_1)+\frac{f'(u)-f'(u_1)}{f'(u_2)-f'(u_1)}(f'(u_2)-f'(u_1))
=f'(u) \;.
\end{align*}
Thus every point on the interpolation $u(x,t)$ satisfies the characteristic
equation \eqref{seibold::eq:characteristic}.\smartqed
\end{proof}
\begin{corollary}[exact solution property]
\label{seibold::cor:exact_solution}
Consider characteristic particles with $x_1(t)<x_2(t)<\dots<x_n(t)$ for
$t\in [t_1,t_2]$. At any time consider the function defined by the interpolation
\eqref{seibold::eq:interpolation_function}. This function is a classical
(i.e.~continuous) solution to the conservation law \eqref{seibold::eq:conslaw}.
In particular, it satisfies the conservation properties given in
Sect.~\ref{seibold::subsec:cons_properties}.
\end{corollary}
\begin{theorem}[TVD]
With the assumptions of Theorem~\ref{seibold::thm:merging_TVD}, the
particle method is total variation diminishing.
\end{theorem}
\begin{proof}
Due to Corollary~\ref{seibold::cor:exact_solution}, the characteristic movement yields
an exact solution, thus the total variation is constant. Particle insertion simply
refines the interpolation, thus preserves the total variation. Due to
Theorem~\ref{seibold::thm:merging_TVD}, merging of particles yields a
new particle with a function value $u_{23}$ between the values of the
removed particles. 
Thus, the total variation is the same as before the merge or smaller.\smartqed
\end{proof}
\begin{remark}
\label{seibold::remark:quadratic_flux}
The method is particularly efficient for quadratic flux functions.
In this case the interpolation \eqref{seibold::eq:interpolation_function} between two
points is a straight line, since $f'$ is linear. Furthermore, the average
\eqref{seibold::eq:nonlinear_average} is the arithmetic mean
$a(u_1,u_2) = \frac{u_1+u_2}{2}$, since $f''$ is constant. Consequently, the function
values for particle insertion and merging can be computed explicitly.
\end{remark}
\begin{remark}
The method has some similarity to \emph{front tracking} by
Holden et al.~\cite{seibold::HoldenHoldenHeghKrohn1988}, and some fundamental
differences. In front tracking, one approximates the flux function by a piecewise
linear, and the solution by a piecewise constant function. Shocks are moved according
to the Rankine-Hugoniot condition. In comparison, our method uses the wave solutions.
Hence, in front tracking everything is a shock; in the particle method, everything is
a wave.
\end{remark}

\section{Entropy}
\label{seibold::sec:entropy}
We have argued in Sect.~\ref{seibold::subsec:interpolation} that due to the
constructed interpolation the particle method naturally distinguishes shocks from
rarefaction fans. In this section, we show that the method in fact satisfies the
entropy condition
\begin{equation}
\eta(u)_t+q(u)\le 0
\label{seibold::eq:entropy_condition}
\end{equation}
if a technical assumption on the resolution of shocks is satisfied. We consider the
Kru\v zkov entropy pair $\eta_k(u) = \abs{u-k}$, $q_k(u) = \sign(u-k)(f(u)-f(k))$.
As argued by Holden et al.~\cite{seibold::HoldenRisebro2002}, if
\eqref{seibold::eq:entropy_condition} is satisfied for $\eta_k$, then it is satisfied
for any convex entropy function.
Relation \eqref{seibold::eq:entropy_condition} implies that the total entropy
$\int\eta_k(u(x))\ud{x}$ does not increase in time for all values of $k$. Using the
interpolation \eqref{seibold::eq:interpolation_function} we show that the numerical
solution obtained by the particle method satisfies this condition. 
\begin{lemma}[entropy for merging]
\label{seibold::lem:entropy_merging}
Consider four particles located at $x_1<x_2=x_3<x_4$, with the middle two to be merged.
We consider the case \mbox{$f''>0$}, i.e.~\mbox{$u_2>u_3$}
WLOG.\footnote{For the case $f''<0$, all following inequality signs must be reversed.}
If the resulting value $u_{23}$ satisfies \mbox{$u_1\ge u_{23}\ge u_4$}, then the
Kru\v zkov entropy does not increase due to the merge.
\end{lemma}
\begin{proof}
We consider the segment $[x_1,x_4]$. Let $u(x)$ and $\tilde u(x)$ denote the
interpolating function before resp.~after the merge. The area under the function is
preserved. We present the proof for $k\le u_{23}$. For $k\ge u_{23}$ the proof is
similar. The interpolating function $u$ is monotone in the value of its endpoints,
thus $u(x)\le\tilde u(x)$ for $x\in [x_2,x_4]$. Since $\abs{x}=x-2\varTheta(-x)$,
where $\varTheta(x)$ is the Heaviside step function, we can write
\begin{align*}
\int_{x_1}^{x_4}\! \abs{u\!-\!k}\ud{x}
&= \int_{x_1}^{x_4}\!\prn{u\!-\!k}\ud{x}
 -2\int_{x_1}^{x_4}\!(u\!-\!k)\varTheta(k\!-\!u)\ud{x} \\
&= \int_{x_1}^{x_4}\!\prn{\tilde u\!-\!k}\ud{x}
 -2\int_{x_2}^{x_4}(u\!-\!k)\varTheta(k\!-\!u)\ud{x} \\
&\ge\int_{x_1}^{x_4}\!\prn{\tilde u\!-\!k}\ud{x}
 -2\int_{x_2}^{x_4}(\tilde u\!-\!k)\varTheta(k\!-\!u)\ud{x}
\ge\int_{x_1}^{x_4}\!\abs{\tilde u\!-\!k}\ud{x} \;.
\end{align*}
Thus, the entropy does not increase due to the merge.\smartqed
\end{proof}
The assumption of Lemma~\ref{seibold::lem:entropy_merging} implies that shocks must be
reasonably well resolved before the points defining it are merged. It is satisfied if
left and right of a shock points are not too far away. In the method, it can be ensured
by an \emph{entropy fix}: A merge is rejected \emph{a posteriori} if the resolution
condition is not satisfied. Then, points are inserted near the shock, and the merge is
re-attempted. It remains to show in future work that with this procedure
Theorem~\ref{seibold::thm:arbitrary_times} still holds.
\begin{theorem}
\label{seibold::thm:entropy}
The presented particle method yields entropy solutions.
\end{theorem}
\begin{proof}
During the characteristic movement of the points the entropy  is constant, since due to
Corollary~\ref{seibold::cor:exact_solution} the interpolation is a classical solution to
the conservation law. Particle insertion does not change the interpolation, thus it 
does not change the entropy. Merging does not increase the entropy if
the conditions of Lemma~\ref{seibold::lem:entropy_merging} are
satisfied.\smartqed
\end{proof}

\section{Numerical Results}
\label{seibold::sec:numerical_results}
\begin{figure}
\centering
\begin{minipage}[t]{.32\textwidth}
\centering
\includegraphics[width=0.99\textwidth]{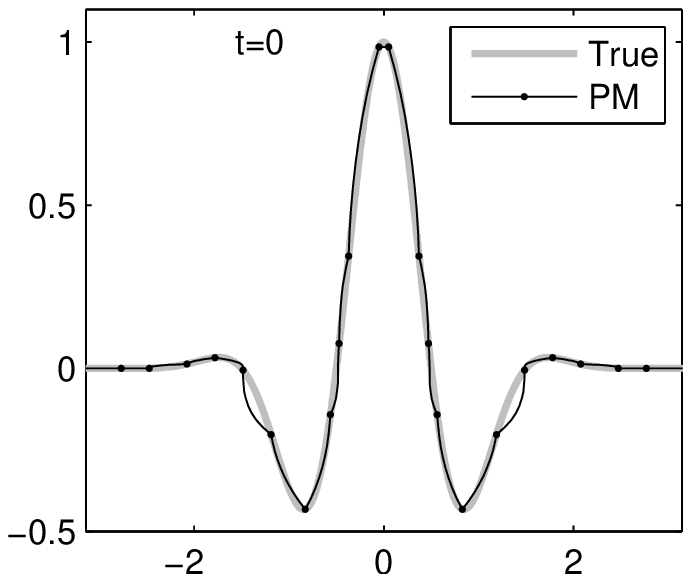}
\end{minipage}
\hfill
\begin{minipage}[t]{.32\textwidth}
\centering
\includegraphics[width=0.99\textwidth]{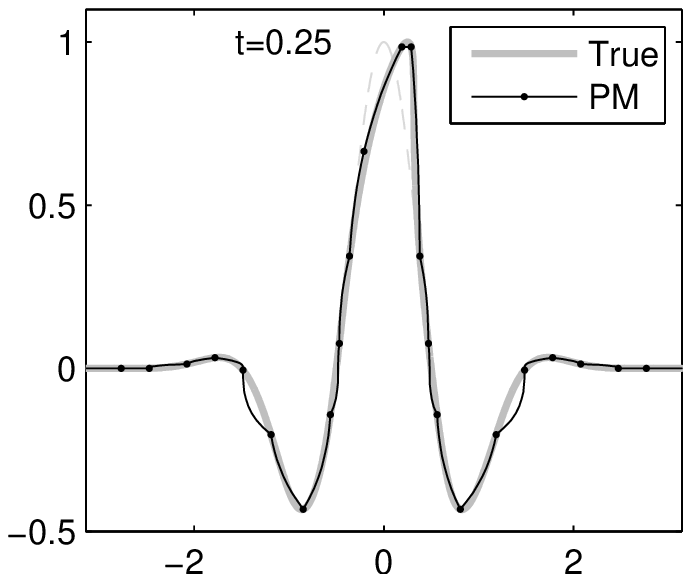}
\end{minipage}
\hfill
\begin{minipage}[t]{.32\textwidth}
\centering
\includegraphics[width=0.99\textwidth]{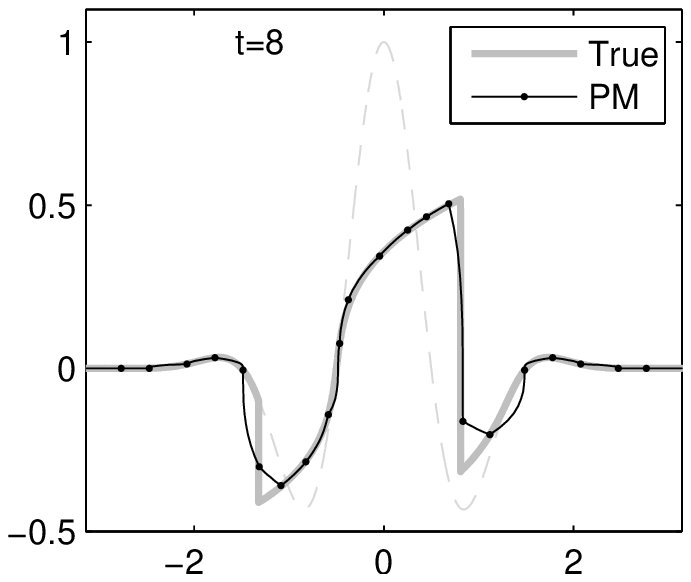}
\end{minipage}
\caption{The particle method for $f(u)=\tfrac{1}{4}u^4$ before and after a shock}
\label{seibold::fig:u4}
\end{figure}

The particle method is particularly well suited for initial conditions that are composed
of similarity solutions. By construction, the movement of the particles yields the exact
solution as long as the solution is smooth. General initial conditions can be approximated
by the interpolation \eqref{seibold::eq:interpolation_function}. Good strategies of
sampling initial conditions shall be addressed in future work.
Figure~\ref{seibold::fig:u4} shows a smooth initial function
$u_0(x) = \exp\prn{-x^2}\cos\prn{\pi x}$ and its time evolution under the flux function
$f(u)=\frac{1}{4}u^4$. The curved shape of the interpolation is due to the nonlinearity
in $f'$. At time $t=0.25$ the solution (obtained by CLAWPACK using 80000 points) is still
smooth, and thus represented exactly on the particles. At time $t=8$ shocks and
rarefactions have occurred and interacted. Although the numerical solution uses only a
few points, it represents the true solution well.

\begin{figure}
\centering
\begin{minipage}[t]{.48\textwidth}
\centering
\includegraphics[width=.99\textwidth]{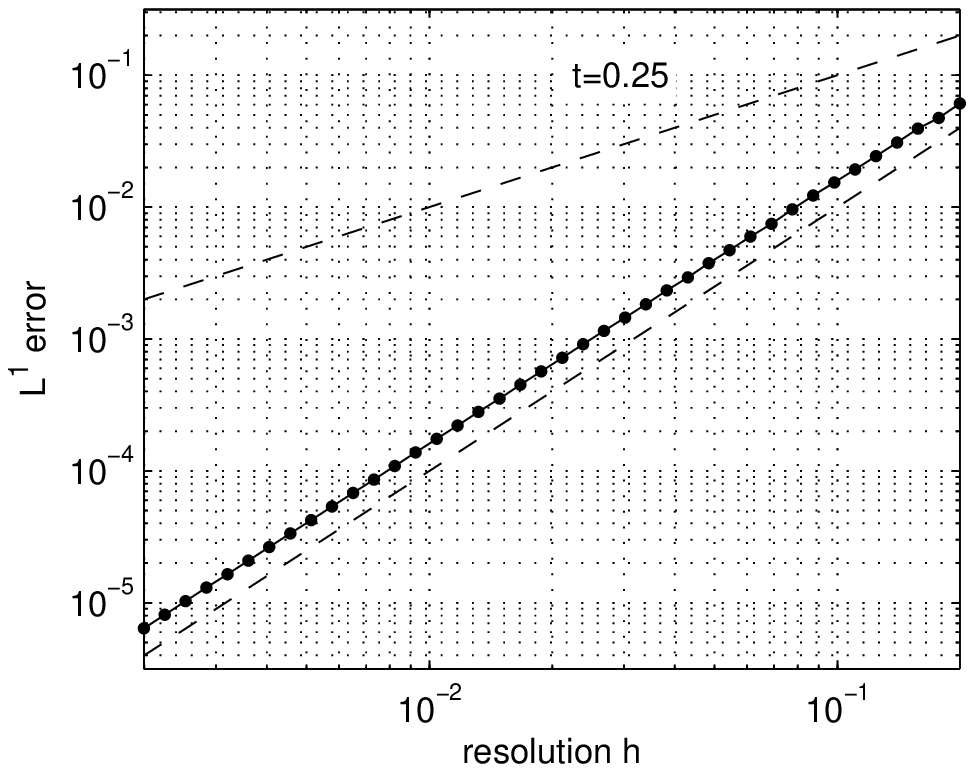}
\end{minipage}
\hfill
\begin{minipage}[t]{.48\textwidth}
\centering
\includegraphics[width=.99\textwidth]{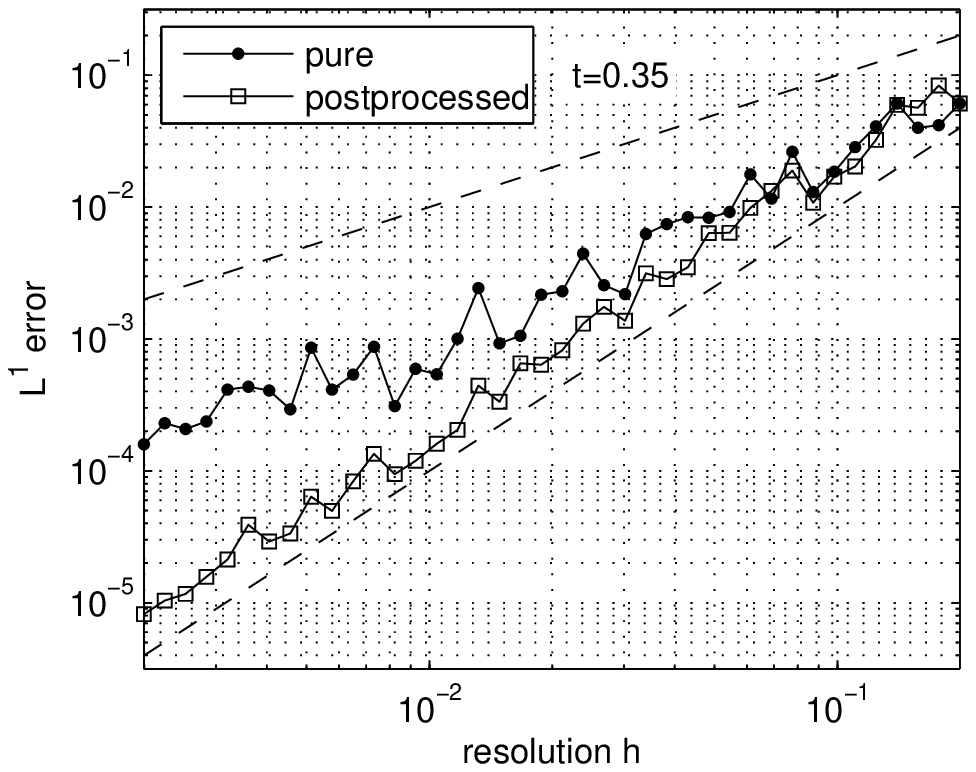}
\end{minipage}
\caption{Error to the correct solution before and after a shock}
\label{seibold::fig:error}
\end{figure}

The accuracy of the particle method is measured numerically. We consider the flux
function and initial conditions as used in Fig.~\ref{seibold::fig:u4}. For a sequence
of resolutions $h$, the initial data are sampled, and the particle method is applied
($d_{\text max}=1.9h$).
Figure~\ref{seibold::fig:error} shows the $L^1$-error to the correct solution
(obtained by a computation with much higher resolution, verified with
CLAWPACK). While the solution is smooth ($t=0.25$), the method is second order accurate,
as is sampling the initial data. After a shock has occurred ($t=0.35$),
the approximate solution (dots) becomes only first order accurate, since the shock has
just been treated by particle management, thus an error of the order
height$\times$width of the shock is made. A postprocessing step (squares) can recover
the second order accuracy: At merged particles, discontinuities are placed so that the
total area is preserved.

\section{Non-Convex Flux Functions}
\label{seibold::sec:inflection_points}
So far we have only considered flux functions with no inflection points (i.e.~$f''$ has
always the same sign) on the region of interest. In this section, we generalize to flux
functions for which $f''$ has a finite number of zero crossings \mbox{$u^*_1<\dots<u^*_k$}.
Between two such points \mbox{$u\in [u^*_i,u^*_{i+1}]$} the flux function is either convex
or concave. We impose the following requirement for any set of particles: Between any two
neighboring particles for which $f''$ has opposite sign, there must be an
\emph{inflection particle} $(x,u_i^*)$. Thus, between two neighboring particles, $f$ has
never an inflection point, and most results from the previous sections apply.
The interpolation between any two particles is uniquely defined
by \eqref{seibold::eq:interpolation_function}. It has infinite slope at the
inflection points, but this is harmless. The characteristic movement of particles is the
same as for flux functions without inflection points. The only complication is merging of
particles when an inflection particle is involved: The standard approach, as presented in
Sect.~\ref{seibold::subsec:particle_management}, removes two colliding points and
replaces them with a point of a different function value. If an inflection particle is
involved in a collision, we must merge points in a different way so that an inflection
particle remains.

\begin{figure}
\centering
\begin{minipage}[t]{.32\textwidth}
\centering
\includegraphics[width=0.99\textwidth]{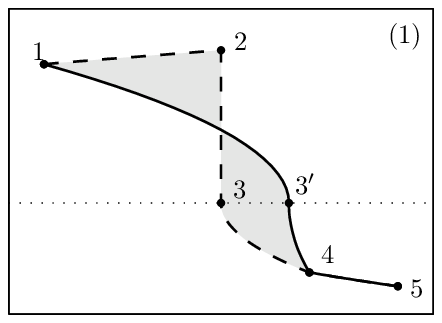}
\end{minipage}
\hfill
\begin{minipage}[t]{.32\textwidth}
\centering
\includegraphics[width=0.99\textwidth]{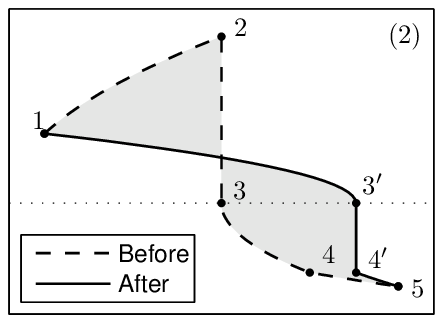}
\end{minipage}
\hfill
\begin{minipage}[t]{.32\textwidth}
\centering
\includegraphics[width=0.99\textwidth]{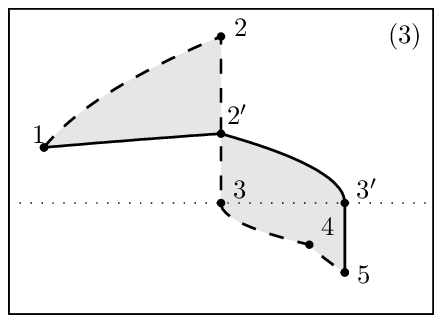}
\end{minipage}
\caption{Particle management around an inflection particle ($f''(u_3)=0$)}
\label{seibold::fig:particle_management_inflection}
\end{figure}

We present one such special merge for dealing with a single inflection point (we do not
consider here the interaction of two inflection points). Also, we only consider a
collision where the positions are exactly the same. Since the inflection particle must
remain (although its position may change), we consider five neighboring particles and
not four as before. Let $(x_i,u_i), i=1,\ldots,5$ be these particle so that $x_2=x_3$,
$f''(u_3)=0$, and (WLOG) $f'''>0$, i.e.~the inflection particle is the slowest. The
other cases are simple symmetries of this situation. We present three successive steps
to finding the final configuration of the particles. Each next step is attempted if the
previous one failed.
\begin{enumerate}
\item
Remove particle 2 and increase $x_3$ to satisfy the area condition. This fails if $x_3$
needs to be increased beyond $x_4$.
\item
Remove particle 2, set \mbox{$x_3=x_4$} and increase both to satisfy the area condition.
This fails if $x_3, x_4$ need to be increased beyond $x_5$.
\item
Remove particle 4, set \mbox{$x_3=x_5$} and find $u_2$ to satisfy the area condition.
This cannot fail if the previous two have failed.
\end{enumerate}
Formally, the resulting configuration could require another, immediate, merge
(since \mbox{$x_3=x_4$} or \mbox{$x_3=x_5$}). However, we need not merge these points
as they move away from each other. The five point particle management guarantees that
in each merging step one particle is removed, thus
Theorem~\ref{seibold::thm:arbitrary_times} holds.

\begin{figure}
\centering
\begin{minipage}[t]{.32\textwidth}
\centering
\includegraphics[width=1.02\textwidth]{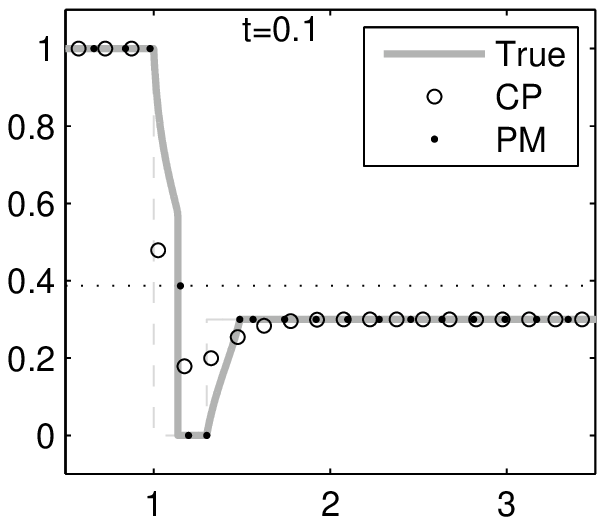}
\end{minipage}
\hfill
\begin{minipage}[t]{.32\textwidth}
\centering
\includegraphics[width=1.02\textwidth]{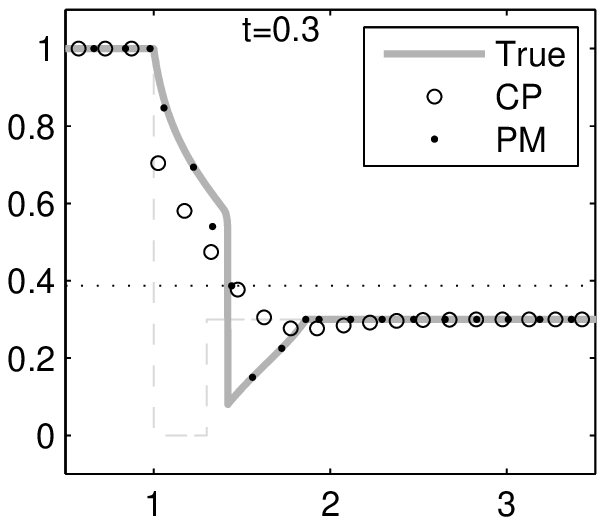}
\end{minipage}
\hfill
\begin{minipage}[t]{.32\textwidth}
\centering
\includegraphics[width=1.02\textwidth]{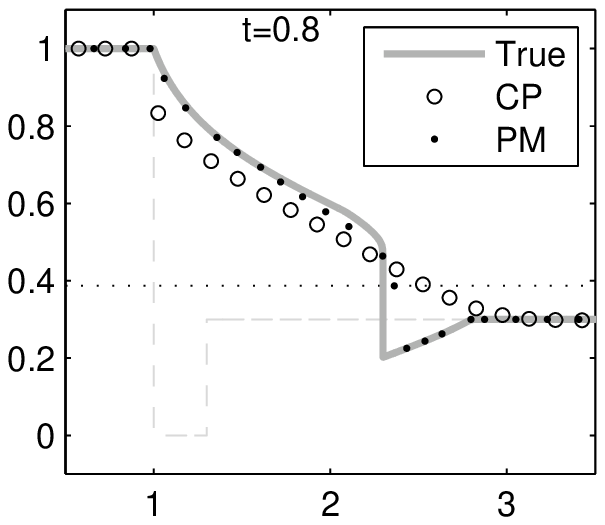}
\end{minipage}
\caption{Numerical results for the Buckley-Leverett equation}
\label{seibold::fig:results_inflection}
\end{figure}

As numerical evidence of the performance, we apply our method to the Buckley-Leverett
equation (see LeVeque \cite{seibold::LeVeque2002}), defined by the flux function
$f(u) = \frac{u^2}{u^2+\frac{1}{2}(1-u)^2}$. It is a simple model for two-phase fluid
flow in a porous medium. We consider piecewise constant initial data with a large
downward jump crossing the inflection point, and a small upward jump. The large jump
develops a shock at the bottom and a rarefaction at the top, the small jump is a pure
rarefaction. Around $t=0.2$, the two similarity solutions interact, thus lowering the
separation point between shock and rarefaction.
Figure~\ref{seibold::fig:results_inflection} shows the reference solution (solid line,
by CLAWPACK using 80000 points). The solution obtained by our
particle method (dots) is compared to a second order CLAWPACK solution (circles) of
about the same resolution. While the finite volume scheme loses the downward jump very
quickly, the particle method captures the behavior of the solution almost precisely.
Only directly near the shock inaccuracies are visible, which are due to the crude
resolution. The solution away from the shock is nearly unaffected by the error at the
shock. Note that although we impose a special treatment only at the inflection
point, the switching point between shock and rarefaction is identified correctly.

\section{Conclusions and Outlook}
\label{seibold::sec:outlook}
We have presented a particle method for 1D scalar conservation laws, which is based on
the characteristic equations. The associated interpolation yields an analytical solution
wherever the solution is smooth. Particle management resolves the interaction of
characteristics locally while conserving area. Thus, shocks are resolved without creating
any numerical dissipation away from shocks. The method is TVD and entropy decreasing,
and shows second-order accuracy. It deals well with non-convex flux functions, as the
results for the Buckley-Leverett equation show.
The particle method serves as a good alternative to fixed grid methods whenever
1D scalar conservation laws have to be solved with few degrees of freedom, but exact
conservation and sharp shocks are desired. An application, which we plan to investigate
in future work, is nonlinear flow in networks (e.g.~traffic flow in road networks). 
For large networks, only a few number of unknowns can be devoted to the numerical
solution on each edge.
We regard the current work as a first step towards a more general particle method. 
Future work will focus on three main directions:
\begin{itemize}
\item\textbf{Source terms:}
Source terms in an equation $u_t+\prn{f(u)}_x = g(x,u)$
could be handled using a fractional step method: 
In each time step, we would first move the particles according to
$u_t+\prn{f(u)}_x = 0$ (including particle management), then change their values
according to an integral formulation of $u_t = g(x,u)$. In the latter step, the
constructed interpolation can be used.
\item\textbf{Systems of conservation laws:}
The particle method is based on similarity solutions of the conservation law. 
For simple systems, such as the shallow water equations in 1D, the analytical solutions
to Riemann problems are known. Two complications arise in the generalization of the
method:
\begin{itemize}
\item
To connect two general states in a hyperbolic system, intermediate states have to be
included.
\item
For a general system it is not clear at which velocity to move the particles.
\end{itemize}
\item\textbf{Higher space dimensions:}
Scalar conservation laws in higher space dimensions can be reduced to 1D problems to be
solved by fractional steps. In principle, this dimensional splitting can be used with
the particle method. However, remeshing would be required between the different spacial
directions, thus the benefits of the meshfree approach would be lost. For the
generalization to a fundamentally meshfree approach in higher space dimensions, the
following problem has to be overcome: In 1D one is never truly meshfree, since the
points have a natural ordering. The method uses this in the interpolation and to detect
shocks. In 2D/3D shocks can occur without particles colliding, as they can move past
each other. Other meshfree methods, such as FPM applied to the Euler equations,
circumvent this issue by using the pressure to regulate shocks.
\end{itemize}

\section*{Acknowledgments}
We would like to thank R.~LeVeque for helpful comments and suggestions. 

\bibliographystyle{amsplain}
\bibliography{seibold_references.bib}



\end{document}